# NEIGHBOR SELECTION AND HITTING PROBABILITY IN SMALL-WORLD GRAPHS

By Oskar Sandberg

*Chalmers University of Technology and Göteborg University*

Small-world graphs, which combine randomized and structured elements, are seen as prevalent in nature. Jon Kleinberg showed that in some graphs of this type it is possible to route, or navigate, between vertices in few steps even with very little knowledge of the graph itself.

In an attempt to understand how such graphs arise we introduce a different criterion for graphs to be navigable in this sense, relating the neighbor selection of a vertex to the hitting probability of routed walks. In several models starting from both discrete and continuous settings, this can be shown to lead to graphs with the desired properties. It also leads directly to an evolutionary model for the creation of similar graphs by the stepwise rewiring of the edges, and we conjecture, supported by simulations, that these too are navigable.

## 1. Introduction.

1.1. *Shortcut graphs.* Starting with the small-world model of Watts and Strogatz [22], rewired graphs have been the subject of much interest. Such graphs are constructed by taking a fixed graph, and randomly rewiring some portion of the edges. Later models of partially random graphs have been created by taking a fixed base graph, and adding "long-range" edges between randomly selected vertices (see [19, 20]). The "small-world phenomenon," in this context, is that graphs with a high diameter (such as a simple lattice) attain a very low diameter with the addition of relatively few random edges.

Jon Kleinberg [11] studied such graphs, primarily ones starting from a two-dimensional lattice, from an algorithmic perspective. He allowed for $O(n)$ long-range edges and found that, not only would this lead to a small diameter, but also that if the probability of two vertices having a long-range edge between them had the correct relation to the distance between them in the grid, the *greedy routing* path-length between vertices was small as well. Greedy routing means, as the name implies, starting from one vertex and searching for another by always stepping to the neighbor that is closest to the destination. That the base graph is connected means that a nonoverlapping greedy path always exists, so the question regards the utility of the long-range contacts in shortening this path. Graphs where one can quickly route between two points using only local information at each step, as with greedy routing, are referred to as *navigable*.

Initially, we will stay in the one-dimensional translation-invariant environment (i.e., with the vertices arranged on a circle). Later sections extend some of the results to other classes of graphs. In general, we will call graphs of the type discussed *shortcut graphs* and use the shorter term *shortcut* for the long-range edges.

1.2. *Contribution.* While Kleinberg's results are important and have been a catalyst for much study, it is not fully understood how the rather arbitrary and precise threshold on the shortcut distribution might arise in practice. In this work, we present an alternative distributional requirement that associates the shortcut distribution with the hitting probabilities of walks under greedy routing. We study this in the canonical case of a single loop, and in a wider setting of graphs induced by the Voronoi tessellations of a Poisson process. We show that distributions that meet a certain criterion which we call being "balanced" have $O(\log^2 n)$ mean routing times, similar to the critical case in Kleinberg's model.

The relationship in this criterion naturally leads to a stepwise rewiring algorithm for shortcut graphs. The Markov chain on the set of possible shortcut configurations defined by this algorithm can easily be seen to have a stationary distribution with balanced marginals. Our analytic results cannot be applied directly in this case, because the stationary distribution has dependencies between the shortcuts at nearby vertices. However, we argue through heuristics and simulation that these dependencies in fact work in our favor, and that graphs generated by our algorithm can be efficiently navigated.

1.3. *Previous work.* The roots of the recent work on navigable graphs are the papers by Jon Kleinberg [10, 11]. Further exposition is given in [1, 16, 17]. Continuum models similar to the ones discussed below have been introduced in [5, 9], and, in a more practical context, [15].

A very different algorithm that appears to produce navigable graphs has been independently proposed in [4], where it is tested by simulation. In [7] the emergence of navigable graphs is discussed in terms of a method for small-world construction without requiring an understanding of the geography, but the method developed is complicated and unnatural. An algorithm similar to that proposed below is present in Freenet [2, 3, 23]—the work below was in part inspired by attempts to place Freenet's algorithms in environments more conducive to analysis.

A recent survey of the field by Kleinberg is [13]. In the final section, he identifies the question of how small-world graphs arise as one of the central questions in the field.

## 2. Preliminaries.

2.1. *Decentralized routing.* The central problem in this area of research is that of *routing* through a graph with only limited knowledge of the graph itself. That is, given two vertices $x$ and $y$ in a (di)graph $G$, we want to find a path connecting $x$ and $y$. In general, the combinatorial problems of finding such a path, and finding the shortest such path, are well-understood problems involving $\Theta(n)$ and $\Theta(n^2)$ steps, respectively. The question becomes more interesting if we allow some (but not all) of the information about the graph to be known when determining the path. In particular, we know the distances between vertices as given by a function $d(x, y)$. With such a distance function, one may define a *decentralized algorithm* (following Kleinberg [11]) as an algorithm which, in each step, uses only information about the distances between vertices already seen in the route and the destination to decide where to go next.

DEFINITION 2.1. A *decentralized algorithm* for finding a path from a point $y$ to $z$ in a graph $G$ associated with a distance function $d: V \times V \mapsto \mathbb{R}$ is defined as follows:

- Let the $S_0 = \{y\}$.
- In step $k$, the algorithm chooses exactly one point in $N(S_{k-1})$ (the set of all neighbors in $G$ of points in $S_{k-1}$) and appends this point to create $S_k$. The choice of $x$ is a possibly random function of the subgraph of $G$ induced by $S_{k-1}$, as well as the distance of all the vertices in $N(S_{k-1})$ to each other and to $z$ as given by $d$. In particular, it may not depend on the rest of $G$.
- The algorithm terminates in step $k$ if $z \in S_k$.

The definition is inspired by the small-world experiments [18] where people were enlisted to forward a letter to a stranger through friend-to-friend links. The people in the experiment knew something about the final recipient (typically where he lived and his occupation), so they could compare how "close" each acquaintance they considered sending the letter to was to the target, but they had no global knowledge of the social network itself.

For a decentralized algorithm to be able to perform better than a random walk, it is necessary that $d(x,y)$ contains some information about the structure of the graph. The extreme of this is where $d(x,y)$ is the graph distance implied by $G$, the minimal distance from $x$ to $y$ in $G$, which we denote $d_G(x,y)$. In this case routing is trivial: proceeding in each step to the neighbor which is closest to $z$ will always find a minimal path. A more typical situation is that $d(x,y)$ gives some, but not complete, information regarding where to go. In particular, we shall say that $d(x,y)$ is *adapted for routing* in a graph $G$, if for any $z$ and $x \in V$, $x$ has a neighbor $y$ such that $d(y,z) < d(x,z)$. When such a distance measure exists, we can route to any point by always choosing such a $y$ as the next step, though the path thus found may be far from optimal.

The common situation is to let $H$ be a fixed graph, and let $G$ be created by randomly augmenting $H$ with further edges in order to create a semirandom graph. It is then trivially true that $d_H(x,y)$ is adapted for routing in $G$. The random edges need not be uniformly distributed, and indeed all the interesting cases arise when the probability of an edge being added between $x$ and $y$ depends on $d_H(x,y)$. Some independence is usually assumed, however, so that the edges previously seen in a route are independent of those in the future. We let $\ell(x,y)$ denote the probability of adding an edge from $x$ to $y$.

Given such a random augmentation of edges, the question arises whether a decentralized algorithm can be found which efficiently routes through a family of graphs. In particular, for a family for finite graphs of bounded degree that are indexed by size, is there a decentralized algorithm such that the expected number of steps of a route between two points is asymptotically small (by which we typically mean that it grows at most poly-logarithmically with the size).

In Kleinberg's original work [11], the underlying graph was $\mathbb{Z}_n^2$ (the family of finite two-dimensional grids) with edges between adjacent vertices, making the $d(x,y)$ the $l^1$ metric (Manhattan distance). He proved that poly-logarithmic routing was possible if $\ell(x,z) = 1/(h_{n,\alpha} x^\alpha)$ with $\alpha = 2$ ($h_{n,\alpha}$ is the distribution's normalizer), but impossible for all other values of $\alpha$. Kleinberg's results also cover the same situation in $\mathbb{Z}_n^d$, in which case the single good value of $\alpha$ is exactly $d$. Similar analysis has been done by others; see, for example, Barriere et al. [1] for thorough analysis of the directed loop, and Duchon et al. [6] for a wider class of graph families. In all these cases (as well as in [12, 14, 15, 21]) it is found that efficient routing is possible when

$$\ell(x,y) \propto \frac{1}{\text{Vol}(B_x(d(x,y)))} \tag{1}$$

where $B_x(r) = \{z : d(x,z) \leq r\}$, or some slight variation thereof. [We will use this notation for the ball, as well as $S_x(r)$ for its boundary throughout the paper.]

Similarly, it turns out that in all these cases, the decentralized algorithm necessary is simply *greedy routing*, which means choosing in each step the unexplored neighbor of the previously explored vertices which is closest to the destination. When $d(x,y)$ is adapted for routing, greedy routing strictly approaches the target with each step and is always successful. The nature of the greedy paths through augmented graphs is the main emphasis of this paper.

The following is a very coarse, obvious, upper bound on the routing time:

OBSERVATION 2.2. *If a distance function $d: V \times V \mapsto \mathbb{R}$ is adapted for routing in a graph in $G = (V, E)$ then greedy routing from $x$ to $z$ takes a number of steps which is at most the cardinality of $\{v \in V : d(v, z) < d(x, z)\}$.*

2.2. *Distribution and hitting probability.* Consider an underlying graph $H = (V, E)$, which may be directed but must be connected in the sense that it contains a direction-respecting path from any vertex to any other. Let $d(x, y)$ be the distance function implied by $H$, and let a random graph $G$ be constructed by augmenting $H$ with one random directed edge starting at each vertex. The edges added by the augmentation will be denoted as $\gamma: V \mapsto V$. We call $\gamma$ a *shortcut configuration*, and let $\Gamma = V \mapsto V$ be the set of all such functions. The general probability space over which we will work is $\Gamma \times V \times V$, where the two copies of $V$ represent the possible starts and destinations of walks. Let $\mathbf{P}$ be a probability measure on this set where the start and destination are chosen uniformly and independently of each other and the configuration is chosen by some *shortcut distribution* $\ell(\gamma)$ which in the independent selection case may be factored into the form $\prod_{x \in V} \ell(x, \gamma(x))$.

On this space, we define $X_Z^Y(t)$ for $t = 0, 1, \ldots$, as a greedy walk in the graph from a uniformly chosen starting point $Y = X_Z^Y(0)$ with a uniformly chosen destination $Z$. To make the greedy walk well defined, we dictate that ties are broken randomly (i.e., if the $m$ closest neighbors to the destination are equally far from it, one is selected uniformly at random). Below, we will in particular be interested in the hitting probability of greedy walks with specific destinations. We define this formally as

(2) $$h(x, z) = \mathbf{P}(X_Z^Y(t) = x \text{ for some } t | Z = z).$$

If $H$ is a translation-invariant graph, then $h(x, z) = h(x - z, 0)$ for some distinguished vertex 0. Thus we will, without further loss of generality, consider the hitting probability as a function of one variable and write $h(x) = h(x, 0)$. Further, we will restrict our analysis to cases where $\ell(x, y)$ and $h(x, y)$ are functions of $d(x, y)$ only (we call this *distance invariance*).

Our results concern relating $h(x)$ to the occurrence of shortcuts between vertices. Immediately, however, we can see that $h(x)$ gives us the expected length of a greedy path. Since such a path can hit each point only once, it follows that if $T$ is the length of a greedy path from a random point to zero, then

$$T = \sum_{x \in V} \chi_{\{X_0^Y(t) = x \text{ for some } t\}}$$

whence it follows that

(3) $$\mathbf{E}[T] = \sum_{x \in V} h(x).$$

We will denote the expected greedy walk length $\tau = \mathbf{E}[T]$.

**3. Rewiring by destination sampling.** Before proceeding to analyze our main model, we present the rewiring algorithm which motivates it. Running the algorithm modifies, in each step, the destinations of some shortcut edges. It is a steady-state algorithm in the sense that the number of edges never changes: it simply shifts the destinations of the single existing shortcut at each vertex.

In the sense that we propose a generative process which might explain why navigable graphs arise, this is similar to the celebrated preferential attachment model for power law graphs of Barabási and Albert. However, it is not a growth model for the graph since the number of vertices and edges never changes, and is thus more similar to the variant of preferential attachment discussed in [8].

The proposed algorithm, which we call *destination sampling*, is as follows:

ALGORITHM 3.1. For a given graph $H = (V, E)$, let $\gamma_s$ be a shortcut configuration at time $s$. From each vertex there is exactly one shortcut. Let $0 < p < 1$. Then $\gamma_{s+1}$ is defined as follows:

1. Choose $y_{s+1}$ and $z_{s+1}$ uniformly from $V$.
2. If $y_{s+1} \neq z_{s+1}$, do a greedy walk from $y_s$ to $z_s$ using $H$ and the shortcuts of $\gamma_s$. Let $x_0 = y_{s+1}, x_1, x_2, \ldots, x_t = z_{s+1}$ denote the points of this walk.
3. For each $x_0, x_1, \ldots, x_{t-1}$ independently with probability $p$ replace its current shortcut with one to $z_{s+1}$ [i.e., let $\gamma_{s+1}(x_i) = z_{s+1}$].

After a walk is made, $\gamma_{s+1}$ is the same as $\gamma_s$, except that the shortcut from each vertex in walk $s + 1$ is with probability $p$ replaced by an edge to the destination. In this way, the destination of each edge is a sample of the destinations of previous walks passing through it (for a realization, see Figure 1). We strongly believe that updating the shortcuts using this algorithm eventually results in a shortcut graph with greedy path-lengths of $O(\log^2 n)$. Though one can relate the stationary regime of this algorithm to the balanced distributions (see below), a rigorous bound has not been proved.

The value of $p$ is a parameter in the algorithm. It serves to disassociate the shortcut of a vertex with those of its neighbors. For this purpose, the lower the value of $p > 0$ the better, but very small values of $p$ will also lead to slower sampling.

3.1. *Markov chain view.* Each application of Algorithm 3.1 defines the transition of a Markov chain on the set of shortcut configurations, $\Gamma$. The Markov chain in question is defined on a finite (if large) state space. If it is irreducible and aperiodic, it thus converges to a unique stationary distribution.

PROPOSITION 3.2. *The Markov chain $(\gamma_s)_{s \geq 0}$ is irreducible and aperiodic.*

PROOF. Aperiodic: There is a positive probability that $y_s = z_s$ in which case nothing happens at step $s$.

Irreducible: We need to show that there is a positive probability of going from any shortcut configuration to any other in some finite number of steps. This follows directly if there is a positive probability that we can "re-point" the shortcut starting at a vertex $y$ to point at a given target $z$ without

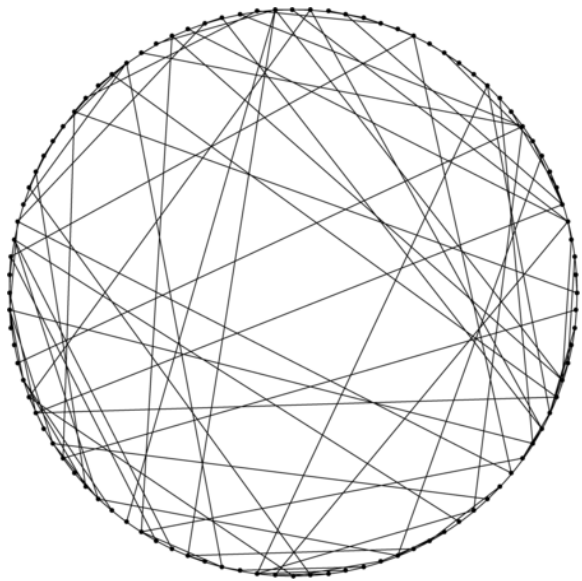

FIG. 1. *A shortcut graph generated by our algorithm ($n = 100$).*

changing the rest of the graph. But the probability of this happening in a single iteration is at least

$$\mathbf{P}(Y=y, Z=z, \text{ and only } y \text{ rewired}) \geq \frac{1}{n}\frac{1}{n}p(1-p)^{n-2} > 0. \qquad \square$$

Thus there does exist a unique stationary shortcut distribution, which assigns some positive probability to every configuration. The goal is to motivate that this distribution leads to short greedy walks.

PROPOSITION 3.3.  *Under the unique stationary distribution of the Markov chain* $(\gamma_s)_{s \geq 0}$

$$\ell(x,z) = \frac{h(x,z)}{\sum_{\xi \neq 0} h(\xi)}.$$

PROOF. As selected by the algorithm, the shortcut from a vertex $x$ at any time is simply a sample of the destination of the previous walks that $x$ has seen. Under the stationary distribution this should not change with time, so

$$\ell(x,z) = \mathbf{P}(Z=z | X_Z^Y(t) = x \text{ for some } t).$$

Using Bayes' theorem, this can be seen as a statement relating $\ell$ to the hitting probability, that is,

$$\ell(x,z) = \mathbf{P}(Z=z | X_Z^Y(t) = x \text{ for some } t)$$
$$= \frac{\mathbf{P}(X_Z^Y(t) = x \text{ for some } t | Z=z)\mathbf{P}(Z=z)}{\sum_{\xi \neq x}\mathbf{P}(X_Z^Y(t) = x \text{ for some } t | Z=\xi)\mathbf{P}(Z=\xi)}.$$

The first multiple in the numerator is the hitting probability $h(x,z)$. The formula then follows from the uniform distribution of $Z$ and translation invariance. $\square$

**4. Balanced shortcut distributions.** We use Proposition 3.3 as the starting point of our analysis, defining the class of all distributions having the same marginal property as follows.

DEFINITION 4.1.  If a graph $H$ with distance function $d(x,y)$ is randomly augmented such that

(4) $$\ell(x,z) = \frac{h(x,z)}{\sum_{\xi \neq 0} h(\xi)}$$

where $h$ is given by (2), then the joint distribution of shortcuts is said to be *balanced*.

We will show for several classes of graphs that this relationship leads to navigability, allowing for a characterization other than that of (1). Besides the relationship with Algorithm 3.1, balance is in some ways a natural requirement. The left-hand side describes the distribution of the destinations of walks that hit the point $x$, so our results simply say that a good choice of shortcuts is one that matches this.

THEOREM 4.2. *For a translation-invariant graph $H$, there exists a balanced distribution which selects shortcuts independently at each vertex.*

PROOF. Like before, we let $\ell(x,y)$ be the marginal probability that $x$ has a shortcut to $y$. The joint distribution is simply the product over all vertices.

For a single walk toward a given $z$, we may view $X_z^Y(t)$ as a Markov chain on the set of vertices, with some transition kernel $P_z(y,x)$. As above, we will set $z=0$ and drop the index in the calculations below without loss of generality. The process hits every point except $z=0$ at most once, and we can let this point be absorbing. The transition kernel $P$ then consists of two mechanisms: either we step to $x$ which is closer to 0 than $y$ because it is the destination of the shortcut from $y$, or we step to one of $y$'s neighbors in $H$ because $y$'s shortcut leads to somewhere from which it is further to 0 than $y$.

Let $N(x)$ be the set of neighbors of $x$ in $H$, and let $L(x) = \{\xi \in V : d(\xi,0) \geq d(x,0), \xi \neq x\}$ be the set of vertices at least as far as $x$ from 0. Also, let $P(x) = \{\xi \in N(x) : d(\xi,0) > d(x,0), (\xi,x) \text{ edge in } H\}$ (the set of "parent" vertices that can greedy route to $x$ in $H$) and $C(x) = \{\xi \in N(x) : d(\xi,0) = d(\xi,0) > d(x,0), (x,\xi) \text{ edge in } H\}$ (the set of "child" vertices that $x$ can route to). Then the transition kernel of the process described is

$$P(y,x) = \begin{cases} 0, & \text{if } d(x,0) \geq d(y,0),\ x \notin C(y), \\ \ell(y,x) + \dfrac{1}{|C(y)|} \sum_{\xi \in L(y)} \ell(\xi), & \text{if } x \in C(y), \\ \ell(y,x), & \text{if } d(x,0) < d(y,0) - 1, \end{cases}$$

for $y \neq 0$. $P(0,x) = \chi_{\{x=0\}}$.

We can thus express the hitting probability for any $x \neq 0$ for a greedy walk as

$$h(x) = \sum_{\xi \in V : \xi \neq x} h(\xi) \ell(\xi,x) + \frac{1}{n-1}$$

$$(5) \qquad = \sum_{\xi : d(\xi,0) > d(x,0)} h(\xi) \ell(\xi,x) + \sum_{\xi \in P(x)} h(\xi) \sum_{\eta \in L(\xi)} \frac{\ell(\eta,\xi)}{|C(\xi)|}$$

$$+ \frac{1}{n-1}.$$

The first two terms in (5) represent the probability that we enter $x$ through either a shortcut or from a parent vertex, respectively, and the last term is the probability that the walk starts at $x$.

Note that, for any $x$, (5) gives a recursive definition of $h(x)$ in terms of the distribution $\ell$. Fix such a distribution $\ell'$. From this we can thus calculate the hitting probabilities $h'(x)$, and define

$$\ell''(x) = \frac{h'(x)}{\sum_{x \in V \setminus \{0\}} h'(x)}.$$

The mapping $\ell' \mapsto \ell''$ is continuous since $\sum_{x \in V} h'(x) > 1$ and maps the simplex of $(n-1)$-valued distributions into itself. Since the simplex is convex, Brouwer's fix-point theorem gives the existence of at least one fix-point $\ell^*$, which is a balanced distribution. □

**5. The directed cycle.** We let $H$ be the directed cycle on $n$ vertices, which will be numbered 0 through $n-1$ such that the edges are directed downward (modulo $n$). The implied distance function (which is not symmetric) is

$$d(x, y) = \begin{cases} x - y, & \text{if } y \leq x, \\ n - y + x, & \text{otherwise.} \end{cases}$$

This environment is perhaps the most natural one for greedy routing, and has previously been the subject of a thorough analysis by [1]. There exists exactly one point at each distance from 0, and greedy routing means selecting the shortcut if and only if its destination lies between 0 and the current position. Equation (5) here simplifies to

$$h(x) = \sum_{\xi=x+1}^{n-1} h(\xi)\ell(\xi - x) + h(x+1) \sum_{\xi=x+2}^{n-1} \ell(\xi) + \frac{1}{n-1}.$$

To prove our result in this environment, we will need the following lemma:

LEMMA 5.1. *If the shortcut configuration is chosen according to a distance-invariant joint distribution, then $h(x)$ is nonincreasing in $x$.*

PROOF. Let $I \subset \Gamma \times V$ be the event consisting of all configurations and starting points such that a greedy walk for 0 hits the point $x + 1$. Now we shift all the coordinates of this set down by one (modulo $n$), and call the translated set $J$. By the definition and distance invariance

$$h(x+1) = \mathbf{P}(I) = \mathbf{P}(J).$$

However, every element in $J$ corresponds to a starting point and shortcut configuration for which the greedy walk hits $x$. To see this, we pick a starting point $y$ and configuration $\gamma$, such that $(\gamma, y) \in I$. This means that there is an integer $m$ and a path $x_0, x_1, \ldots, x_m$ such that $x_0 = y$, $x_m = x + 1$ and either

$$n - 1 \geq \gamma(x_i) > x_i \quad \text{and} \quad x_{i+1} = x_i - 1$$

or

$$x_i > \gamma(x_i) \geq x + 1 \quad \text{and} \quad x_{i+1} = \gamma(x_i)$$

for all $i = 0, 1, \ldots, m$. The corresponding configuration in $J$ has a similar path $x'_0, \ldots, x'_m$ ($x'_i = x_i - 1$) where $x'_0 = y - 1$, $x'_m = x$ and either

$$n - 2 \geq \gamma(x'_i) > x'_i \quad \text{and} \quad x'_{i+1} = x'_i - 1$$

or

$$x'_i > \gamma(x'_i) \geq x \quad \text{and} \quad x'_{i+1} = \gamma(x_i)'$$

for all $i = 0, 1, \ldots, m$. This means that starting in $y - 1$ will cause the greedy walk to hit $x$. [Note that not every configuration and starting point that cause greedy walks to hit $x$ are necessarily in $J$, since $\gamma(x'_i)$ must be less than $n - 2$ rather than $n - 1$ in the first line.]

It now follows directly that

$$\mathbf{P}(J) \leq h(x). \qquad \square$$

We can now show that greedy routing here has taken a similar number of steps to the critical case in Kleinberg's model.

THEOREM 5.2. *For every $n = 2^k$ with $k \geq 3$, the shortcut graph with shortcuts selected independently according to a balanced distribution has an expected greedy routing time*

$$\tau \leq 2k^2.$$

The proof method is similar to that introduced by Kleinberg for augmentations described by (1) links, but the implicit definition of the shortcut distribution requires a somewhat more involved approach.

PROOF OF THEOREM 5.2. Assume that $\tau > 2k^2$. We will show that for $k$ sufficiently large this always leads to a contradiction.

To start with, divide $\{1, 2, \ldots, n-1\}$ into at most $k$ disjoint *phases*. Each phase is a connected set of points, each successively further from the destination 0, and they are selected so that a greedy walk is expected to spend the same number of steps in each phase. Thus, the first phase is the interval $F_1 = \{1, \ldots, r_1\}$ where $r_1$ is the smallest number such that

$$\ell(F_1) = \sum_{\xi \in F_1} \ell(\xi) \geq 1/k.$$

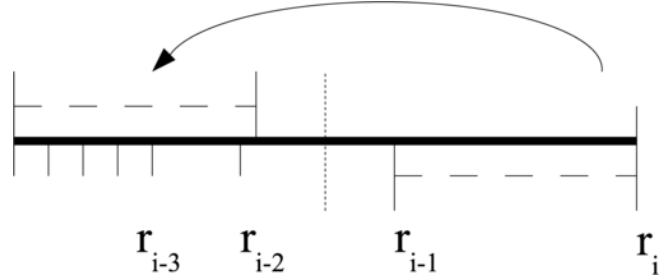

Fig. 2. *Illustration of the proof of Theorem 5.2. If a phase covers less than half of the "remaining ground," then any shortcut of the same distance from $r_i$ as the 0 is from the points in the phase takes us out of the phase.*

The second phase is defined similarly as the shortest interval $\{r_1+1,\ldots,r_2\}$ such that $\ell(F_2) \geq 1/k$. Let $m$ be the total number of such intervals which can be formed, and let $F_R$ denote the remaining interval $\{r_m+1,\ldots,n-1\}$ which could be empty. By construction $\ell(F_R) < 1/k$ and the total number of phases, including $F_R$, is at most $k$.

Before proceeding, we need to bound by how much $\ell$ of the different phases can deviate from one another since this will also tell us by how much the expected number of steps in each phase can differ. From (4) and the assumed lower bound of $\tau$, it follows that

$$\ell(x) = \frac{h(x)}{\tau} \leq \frac{1}{2k^2}$$

for all $x$. This implies that

$$\frac{1}{k} \leq \ell(F_i) \leq \frac{1}{k} + \frac{1}{2k^2}$$

for all $i \in \{1,\ldots,m\}$, and thus

(6)  $$\ell(F_i) \leq \left(1 + \frac{1}{2k}\right)\ell(F_j)$$

for all $i,j \in \{1,\ldots,m\}$. It also gives $m \geq k^2/(k+1) - 1$.

Consider now $F_m = \{r_{m-1}+1, r_{m-1}+2, \ldots, r_m\}$ and let $L = \{0,1,\ldots,r_{m-1}\}$. Our goal is to show that, from any point in $F_m$, there is a considerable probability of having a shortcut to $L$. We know that $r_m \leq n$. Assume that $r_{m-1} \geq r_m/2$. $F_m$ then covers less than half of the distance from $r_m$ to the target. In particular

$$\{r_m - F_m\} \subset L.$$

Thus, if $r_m$ has a shortcut with destination in $\{r_m - F_m\}$, any walk which hits $r_m$ will leave $F_m$ in the next step (see Figure 2). We thus know that

$$\ell(r_m, L) \geq \ell(r_m, \{r_m - F_m\}) = \ell(F_m) \geq 1/k.$$

Lemma 5.1 tells us that the probability of having a shortcut to $L$ cannot decrease for points less than $r_m$, so for each vertex that the walk hits within $F_m$, there is an independent probability of at least $1/k$ of leaving $F_m$ in the next step. This means that the expected number of steps the walk can take in $F_m$ is at most $k$.

The expected number of steps in a phase is $h(F_i) = \tau \ell(F_i)$ so, by (6), it then holds that

(7) $$h(F_i) \leq (1 + 1/2k) h(F_m) \leq k + 1/2$$

for all $i \in \{1, \ldots, m\}$ and also for $F_R$. There are at most $k$ phases, so this implies that $\tau \leq k^2 + k/2$, which contradicts our assumption for all $k \geq 2$.

Thus the original assumption implies that $r_{m-1} \leq r_m/2 \leq n/2$. But by an identical argument for $F_{m-1}$, we can show that $r_{m-2} \leq r_{m-1}/2$. It follows by iteration that

$$r_i \leq \frac{1}{2^{m-i}} n$$

and in particular

$$r_1 \leq \frac{1}{2^{m-1}} n \leq 2^{(k+2)/(k+1)} \leq 4.$$

This means that $F_1$ contains at most four points, so $h(F_1) \leq 4$ and thus, again by (6), $h(F_i) \leq 5$ for all $i$. For $k \geq 3$ this contradicts the original assumption. This completes the proof. $\square$

Theorem 5.2 gives us an alternative distributional criterion for attaining $O(\log^2 n)$ expected greedy path-lengths. Since Kleinberg showed that this cannot hold for many distributions, the balanced distributions must be "close" to the critical, harmonic decay. More specifically, drawing on the proofs that navigability is not possible for most case graphs, we can see that there cannot exist $\delta > 0, \epsilon > 0$ and $N \in \mathbb{N}$ such that $\ell(\{1, 2, \ldots, n^\delta\}) \leq n^{-\epsilon}$ for the cycles of size $n \geq N$. This would be the case if the tails of the distributions dominated a power law $(x^{-\alpha})$ decay with exponent $\alpha < 1$. Similarly, there cannot exist (possibly different) $\delta > 0, \epsilon > 0$ and $N \in \mathbb{N}$ such that $\ell(\{n^\delta, n^\delta + 1, \ldots, n-1\}) \leq n^{-\epsilon}$ for the cycles of size $n \geq N$, as would be the case if the tails were dominated by a power law with exponent $\alpha > 1$.

**6. Delaunay graphs.** The small-world theory is not necessarily limited to situations where vertices are placed in a fixed grid. In this section, we will let the vertices be points of a spatial Poisson process, and the distance function be the Euclidean metric. For simplicity, we will relax our requirements a little and let the graphs have degrees bounded in expectation, rather than uniformly.

Let $S^d$ be the $d$-dimensional surface of a $d+1$ ball with radius such that the volume/area of $S^d$ is 1. We let $V = \{x_i\}$ be the $N$ points of a homogeneous Poisson process with intensity $\lambda = n^d$ in this space. From this Poisson process we may construct the Voronoi tessellation, that is, the collection of cells $C(x_i)$ where

$$C(x_i) = \left\{ y \in \mathbb{R}_s^d : |y - x_i| = \min_{z \in V} |z - x_i| \right\}.$$

$C(x_i)$ is that part of the space which is as close to $x_i$ as to any other point. The Voronoi cells are closed convex polyhedrons that border other cells along each side, thus overlapping on sets of Lebesgue measure zero.

The tessellation induces a graph $G$ with vertices $V$ (known as the Delaunay graph) as follows. Let $G = (V, E)$, where $(x, y) \in E$ if and only if $C(x)$ and $C(y)$ intersect in an infinite number of points (this is a.s. equivalent to intersecting in at least one point). Intuitively, this is the graph that connects a vertex $x$ to all its neighbors in the tessellation. This Delaunay graph is a natural base graph for greedy routing among the points.

LEMMA 6.1. *Let $\{x_i\}$ be any point-set in $S^d$, and $G$ its Delaunay graph. Then the Euclidean metric $d(x, y) = |x - y|$ is adapted for routing in $G$.*

PROOF. We must prove that for all $x \neq z \in V$, there exists $y \in V$ (which may be $z$) such that $(x, y) \in E$ and $|x - z| > |y - z|$. Consider the line $xz$. Let $w$ be the first point we encounter as we move from $x$ along $xz$, satisfying $w \in C(y)$ for some $y \neq x$ ($w$ is well defined since the cells are compact).

It is clear that $w \in C(y)$ for some $y$ such that $x$ and $y$ are connected in the Delaunay graph [$C(x)$ must border at least one cell that it meets at $w$]. Clearly, $|y - w| = |y - w|$ since $w$ is in both cells. Thus

$$|y - z| < |y - w| + |w - z|$$
$$= |x - w| + |w - z| = |x - z|$$

where the strict inequality follows from the fact that $w$ is not on the line $yz$. □

Given this graph, we consider augmentations that allow for fast routing. A direct approach would be to connect a given vertex to any other with a probability depending on the distance between them, but this leads to complications regarding dependencies between the progress made at each step (though not insurmountable ones; see [5] for such an approach in a similar environment).

Instead, we augment the graph as follows. For each vertex $x \in V$, let $\{n_i(x)\}_{i=1}^{N(x)}$ be the points of a nonhomogeneous Poisson process given by the measure $\mu_x(A) = \ell_x(A \setminus C(x))$ for some *shortcut measure* $\ell$ on the Borel sets of $S^d$, and $\ell_x(A) = \ell(A - x)$.

We then augment $G$ by adding a directed edge from $x$ to $y$ if $n_i(x) \in C(y)$ for any $i = 1, \ldots, N$.

LEMMA 6.2. *If $x, z \in V$ and $|z - n_i(x)| \leq |z - x|/4$ for some $i = 1, 2, \ldots, N$, then $x$ has a shortcut $y \in V$ (which may be $z$) such that $|z - y| \leq |z - x|/2$.*

PROOF. Let $w$ be such an $n_i(x)$. With probability 1 it is in exactly one cell $C(y)$. If $y = z$, then $x$ has a shortcut to $z$; otherwise $|w - y| < |w - z|$. In the latter case,
$$|z - y| \leq |z - w| + |w - z| < 2|w - z| \leq |x - z|/2. \qquad \square$$

6.1. *Kleinberg augmentation.* To motivate the model, we first show that augmentation along the lines of Kleinberg's model allows for an $O(\log^2(n))$ bound on the routing time. That is, as in (1), we let the augmentation be given by the measure

$$(8) \qquad \ell(A) = \int_A \frac{dr}{\log n \, \mathrm{Vol}(r)}$$

where $\mathrm{Vol}(r)$ is the volume of a ball of radius $r$ in $S^d$. The measure is defined on sets $A \in S^d \setminus \{0\}$.

Before proving a lower bound on the expected routing type, we need to ensure that we are not adding an unbounded number of edges.

LEMMA 6.3. *The expected number of shortcuts added to each vertex under augmentation with intensity (8) is bounded by a constant.*

PROOF. First note that $\mathbf{E}[\#\text{shortcuts added to } x] \leq \mathbf{E}[N(x)]$. Now, let $R(x) = \inf\{|y - x| : y \in V, y \neq x\}$. If $R(x) = \delta$, then all points within distance $\delta/2$ of $x$ are in $C(x)$. Thus

$$\mathbf{E}[N(x) \mid R(x) = \delta] \leq \frac{1}{\log n} \int_{S^d \setminus B_{\delta/2}(x)} \frac{1}{\mathrm{Vol}(x - y)} \, dy$$

$$\leq \frac{1}{\log n} \int_{\delta/2}^{1} \frac{1}{r} \, dr = \frac{\log(2/\delta)}{\log n}.$$

Hence, and since $\mathbf{E}[N(x) \,|\, R(x) = \delta]$ is decreasing in $\delta$,

$$\begin{aligned}
\mathbf{E}[N(x)] &= \int_0^1 E[N(x)\,|\,R(x)=\delta]f_{R(x)}(\delta)\,d\delta \\
&\leq \mathbf{E}[N(x)\,|\,R(x)=1/n]\mathbf{P}(R(x) \geq 1/n) \\
&\quad + \int_0^{1/n} \frac{\log(2/\delta)}{\log n} n^d S(\delta) e^{-n^d \operatorname{Vol}(\delta)}\,d\delta \\
&\leq 2 + \frac{n^d S(1/n)}{\log n} \int_0^{1/n} \log(2/\delta)\,d\delta \\
&\leq c,
\end{aligned}$$

where $S(\delta)$ is the area of a sphere of radius $\delta$, and $c$ is a constant independent of $n$.  $\square$

The proof of the following theorem uses the by-now standard argument from [11].

THEOREM 6.4. *For every $n$ sufficiently large, the shortcut graph created by augmenting the Poisson–Delaunay graph with intensity (8) has an expected greedy routing time of $O(\log^2 n)$.*

PROOF. Let the route currently be at the vertex $x$, such that $|x - z| = d > 1/n$. Let $B$ be the event that $|n_i(x) - z| \leq d/4$ for some $i$. Then

$$\mathbf{P}(B) \geq \frac{\operatorname{Vol}(d/4)}{\operatorname{Vol}(3d/4)\log n} = \frac{c}{\log n}.$$

By Lemma 6.2, if such a $n_i(x)$ exists, then $x$ has a neighbor within distance $d/2$ of $z$, and greedy routing at least halves the distance to $z$ in the next step. If $B$ fails to occur, then we know by Lemma 6.1 that greedy routing can still progress to a point closer to the destination, and whether or not $B$ occurs is independent of previous steps. Thus the expected number of steps until the distance to the target is halved is $O(\log n)$, which together with Lemma 2.2 proves the result.  $\square$

6.2. *Balanced augmentation.* In order to derive a result similar to Theorem 5.2 for the Delaunay setting, we will need to redefine the "balanced distribution" somewhat. In particular, we need to marginalize over the positions of the Poisson points.

Let the hitting measure of $A \subset S^d \backslash \{0\}$ be defined by

$$h_z(A) = \mathbf{E}[\text{number of } t \text{ s.t. } X_Z(t) \in A \,|\, Z = z]$$

where $X_z(t)$ is the greedy routing process as above, and the existence of a point at $z$ is included in the conditioning. Note that, by the translation invariance of the construction, $h_z(A) = h_0(A - z)$.

We call a distribution *Poisson-balanced* if

(9) $$\ell(A) = \frac{h_0(A)}{\tau}$$

where $\tau = \mathbf{E}[\text{length of a greedy walk}] = h_z(S^d \backslash \{0\})$.

LEMMA 6.5. *There exists a Poisson-balanced distribution.*

PROOF. The proof is similar to the discrete case. A given shortcut measure $\ell'$ induces a hitting measure $h'_0(A)$, which in turn gives rise to a measure $\ell''$ via (9). If we let $L$ be the space of measures of total probability 1 on $S^d \backslash \{0\}$ equipped with the total variation metric

$$d_{\text{TV}}(\mu, \nu) = \sup_{A \in \mathcal{B}(S^d \{0\})} |\mu(A) - \nu(A)|,$$

then the mapping $\ell' \mapsto \ell''$ is a mapping from $L$ to itself. $L$ is convex and compact, so it suffices to show that the mapping is continuous for us to apply Brouwer's fix-point theorem.

Since we know that $\tau > 1$, the second step of the mapping is certainly continuous. The first is also, since the hitting probability depends only on a finite number of random variables with distribution depending on $\ell'$. Formally:

Take $\epsilon > 0$ and any $m = n^d$. Let $\ell_1$ and $\ell_2$ be two shortcut measures. Without loss of generality, we assume that $\ell_2 \geq \ell_1$, and we let $d_{\text{TV}}(\ell_1, \ell_2) \leq \epsilon'$ where

$$\epsilon' \leq \frac{\epsilon}{3m \max((e-1)n, \log(3m/\epsilon))}.$$

We couple the routing processes $X_0^1$ and $X_0^2$ by letting them use the same set of Poisson process distributed vertices $V$, and the same starting point $z$. At each $x \in V$, we construct shortcuts $n_i(x)$ according to $\ell_1$ which both processes may use, and then add an additional set of shortcuts $\{n_i^2(x)\}$ according to $\ell_2 - \ell_1$, which only $X_0^2$ may use.

It follows that for any $x$, the cardinality of $\{n_i^2(x)\}$, $N(x)$, is dominated by a $\text{Poi}(\epsilon')$ random variable, so

$$\mathbf{P}(N(x) = 0) \leq 1 - e^{-\epsilon'} \leq \epsilon'.$$

Let $B$ be the event that a given vertex $x$ in $V$ has $N(x) > 0$. Then

$$\mathbf{P}(B) \leq \mathbf{P}(B \,|\, |V| \leq (e-1)m + q) + \mathbf{P}(|V| \leq (e-1)m + q)$$
$$= ((e-1)m + q)\epsilon' + e^{-q}$$
$$\leq \frac{\epsilon}{m},$$

where the last inequality follows from setting $q = \log(3m/\epsilon)$.

Now let $H_1(A)$ and $H_2(A)$ be the number of points reached in a subset $A \subset S^d \backslash \{0\}$ by the respective processes. If the $h_0^1$ and $h_0^2$ are the respective hitting probabilities, then

$$|h_0^1(A) - h_0^2(A)| = \mathbf{E}[|H_1(A) - H_2(A)|]$$
$$= \mathbf{E}[|V|]\mathbf{P}(H_1(A) \neq H_2(A)) \leq \epsilon$$

since $H_1(A) = H_2(A)$ if no vertex has different shortcuts in the two cases. This completes the proof. □

In order to bound the routing time in this case, we will need the following geometrical fact.

LEMMA 6.6. *There exists $q \in (0,1)$ such that, if $x$ and $y$ are points in a $S^d$, satisfying $(3/4)\delta < |x - y| \leq \delta$, and $(3/4)\delta < r \leq \delta$, then the portion of the sphere $S_r(y)$ which lies inside $B_{(3/8)d}(x)$ is at least $q$. The constant $q$ depends on $d$ but not on $\delta$ and $r$.*

This follows directly from the fact that the statement is independent of scale. In one dimension $q = 1/2$ trivially, and in two it can easily be seen that it is at least $1/8$; see Figure 3.

THEOREM 6.7. *For every sufficiently large $k$ and $n = (\frac{4}{3})^k$, the shortcut graph created by augmenting the Poisson–Delaunay graph with a Poisson-balanced distribution has an expected greedy routing time $\tau \leq \frac{k^2}{q}$.*

PROOF. We let $X_0^Y$ be the routing process for zero, and define $h_0$ on $S^d \backslash \{0\}$ as above. We then divide $S^d \backslash \{0\}$ into $k$ phases of the form $F_i = \{x \in$

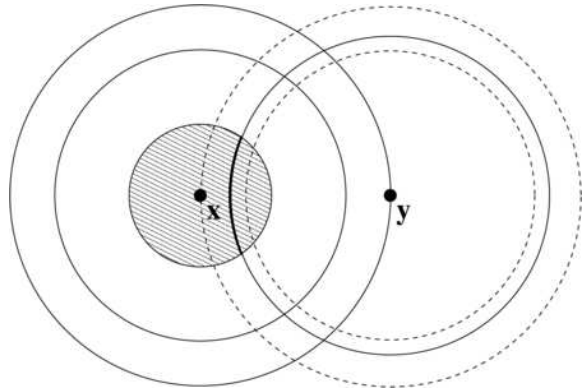

FIG. 3. *The circle around $y$ intersects the ball around $x$ in at least $1/8$ of its points.*

$S^d : r_{i-1} < |x| \leq r_i\}$, where $r_0 = 0$ and each subsequent $r_i$ is defined so that

$$h_0(F_i) = \frac{\tau}{k}.$$

For any phase $F_i$, assume that $r_{i-1} \geq \frac{3}{4} r_i$. Let $x$ be a vertex in $F_i$. A portion $q$ of the area of each spherical "level" in $x + F_i$ lies in $L_i = B_0((3/8)r_i)$ by Lemma 6.6. By rotational invariance it follows that $\ell_x(L_i) = q\ell_x(x + F_i) = q\ell(F_i)$, so if $B$ is the event $x$ has a shortcut destination $n_i(x)$ closer than $r_{i-1}/2$, then

$$\mathbf{P}(B) = \ell_x(B_0((3/8)r_i)) \geq q \frac{h_0(F_i)}{\tau} = \frac{q}{k}.$$

By Lemma 6.2, if such a $n_i(x)$ exists, then $x$ has a neighbor within distance $d/2$ of $z$, and greedy routing at least halves the distance to $z$ in the next step. If $B$ fails to occur, then we know by Lemma 6.1 that greedy routing will progress to a vertex closer to the target. The event $B$ is independent of previous steps. Thus $h_0(F_i) \leq \frac{k}{q}$, whence $\tau \leq \frac{k^2}{q}$.

If, on the other hand, $r_{i-1} \leq \frac{3}{4} r_i$ for all $i$, then

$$r_1 \leq \frac{3^{k-1}}{4} = \frac{4}{3n}.$$

Let $N$ be the number of vertices in $F_1$. By Observation 2.2

$$h_0(F_1) \leq \mathbf{E}[N] = \frac{\text{Vol}(r_1)}{n^d} = c.$$

It follows that $\tau \leq ck$, so the result holds when $k > cq$. □

**7. The rewiring algorithm revisited.** Proposition 3.3 shows that, under the stationary distribution of the destination sampling algorithm introduced above, the marginal shortcut distribution at each point is balanced, and it is thus tempting to apply Theorem 5.2 to bound the greedy path-length. However, that theorem assumed that the shortcuts had been chosen independently at each vertex, which is not the case under Algorithm 3.1 which originally motivated the work. Showing that these dependencies do not negatively affect routing is an open problem, which we discuss in general terms in this section.

There are two sources of dependencies between the shortcuts of neighboring vertices. First, there is a chance that they sampled the destination of the same walk. When $p$ is large, this dependency is substantial, and we see a highly detrimental effect even in the simulations. By using a small $p$, however, this dependence is muted. Another, more subtle dependence has to do with the way the shortcuts of vertices around a vertex $x$ may affect the destinations of the walks $x$ sees. In the directed cycle, if $x + 1$ has a shortcut to $x - 10$, that will make it less likely for $x$ to see walks for places "beyond" $x - 10$, since many such walks will have followed the shortcut at $x + 1$, and thus skipped over $x$.

The first source of dependence, that of sampling from the same walk, can be handled by modifying the algorithm to make sure we do not sample more than once for each walk. Take $p \leq 1/n$ and, once a walk is completed, choose to update exactly one of its links with probability $pw$ where $w$ is the length of the walk. Which link to update is then chosen uniformly from the walk. This way, the probability that a vertex updates its shortcut when hit by a walk is still always $p$, but we never sample two shortcuts from the same walk. The modified algorithm is less natural, but clearly a good approximation of the original for small values of $p$. Although it is more complicated, it is easier to analyze, since it allows for the simplifying assumption that each edge is chosen from a different greedy walk.

The other dependencies are more complicated, and there is no easy way to modify the algorithm to remove them. However, it is worth noting that it is hard to see why these dependencies (unlike the first type) would be destructive for greedy routing. In fact it makes sense that, if $x$ in our example gets few walks destined beyond $x - 10$ because of the shortcut present at $x + 1$, then it should also choose a shortcut to beyond $x - 10$ with a smaller probability.

In the proof of Theorem 5.2 we use independence only to show that if the probability of having a shortcut out of a phase at the very furthest point is $\rho$, then the expected number of steps in the phase is bounded by $1/\rho$. There is little reason to believe this would not hold under the modified algorithm, since if the link from the furthest point does not take us out of the phase, then it either goes to a point within the phase, or overshoots the destination. If it goes to a point within the phase, then we follow it, and the presence of that shortcut should not interfere with those we see in the future. If, on the other hand, it overshoots, then by the argument above it should make it more likely that the succeeding ones do not do so, giving us a better probability of leaving the phase than in the independent case.

Formalizing the requirements on the dependence, and proving that our stationary distribution indeed has the necessary properties, is the main open problem which we have yet to resolve.

7.1. *Computer simulation.* Simulations indicate that the algorithm gives results which scale as desired in the number of greedy steps, and that the distribution approximates $1/(h_{n,d}\, d(x,y))$ for the one-dimensional grid.

The results in the directed one-dimensional grid can be seen in Figure 4. To get these results, the graph is started with no shortcuts, and then the algorithm is run $10n$ times to initialize the references. The value $p = 0.1$ is used. The greedy distance is then measured as the average of 100,000 walks, each updating the graph according to the algorithm. The effect of running the algorithm, rather than freezing one configuration, seems to be to lower the variance of the observed value.

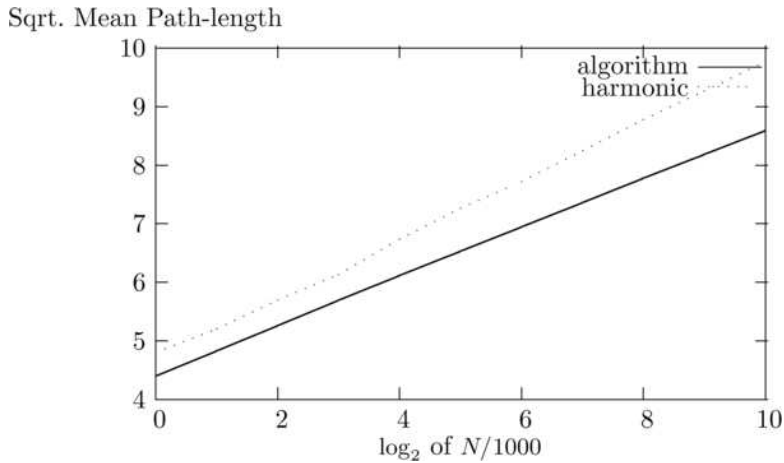

FIG. 4. *The expected greedy walk time of the selection algorithm, compared to selection according to harmonic distances, in a cycle.*

The square root of the mean greedy distance increases linearly as the graph size increases exponentially, just as we would expect. In fact, as can be seen, our algorithm leads to better simulation results than choosing from Kleinberg's distribution. Doubling the graph size is found to increase the square root of the greedy distance by about 0.41 when links are selected using our algorithm, compared to an increase of about 0.51 when Kleinberg's model is used. [In fact, in Kleinberg's model we can use (5) to calculate numerically exact values for $\tau$, allowing us to confirm this figure.]

In Figure 6 the marginal distribution of shortcut lengths is plotted. It is roughly harmonic in shape, except that destination sampling creates fewer links of length close to the size of the graph. This may be part of the reason why it is able to outperform Kleinberg's model: while the latter is asymptotically correct, our algorithm takes into account finite size effects. (This reasoning is similar to that of the authors of [4]. Like them, we have no strong analytic arguments for why this should be the case, which makes it a tenuous argument at best.)

The algorithm has also been simulated to good effect using base graphs of higher dimensions. Figure 5 shows the mean greedy distance for two-dimensional grids of increasing size. Here also, the algorithm creates configurations that seem to display square logarithmic growth, and which perform considerably better than explicit selection according to Kleinberg's model.

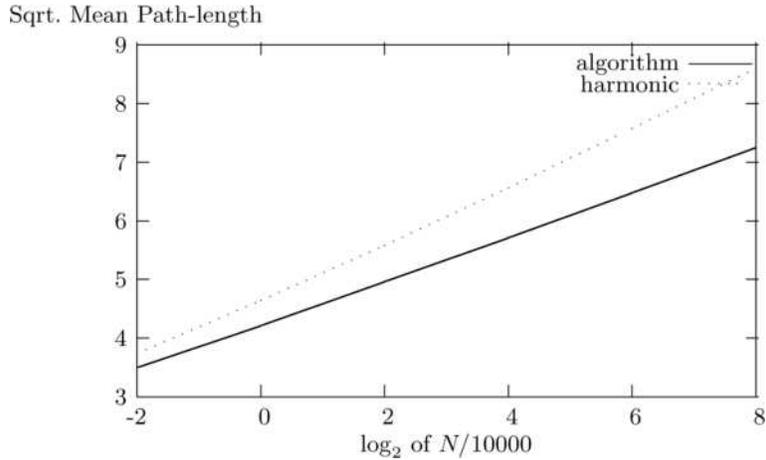

Fig. 5.  *The expected greedy walk time of the selection algorithm, compared to selection according to harmonic distances, in a two-dimensional base grid.*

**8. Conclusion.**  The study of navigable graphs is still in its infancy, but many interesting results have already been found, and the practical relevance to such fields as computer networks is beyond doubt. In this paper we have presented a different way of looking at the dynamics that cause graphs to be navigable, and we have presented an algorithm which may explain how navigable graphs arise naturally. The algorithm's simplicity also means that it can be useful in practice for generating graphs that can easily be searched, an important property for many structures on the Internet.

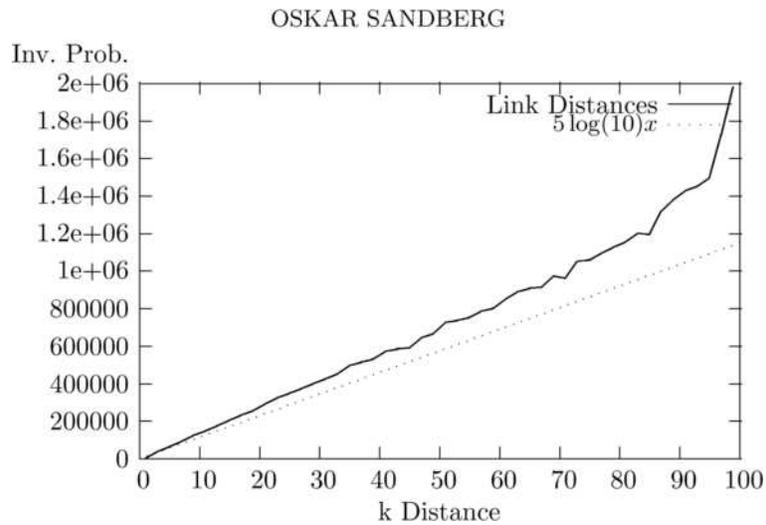

Fig. 6.  *The inverse of the distribution of shortcut distances, with $n = 100\,000$, $p = 0.1$. The straight line is the inverse of the harmonic distribution.*

While many questions about these graphs in general, and our results in particular, remain unanswered, the prospects for going further with this work seem good. We are hopeful that these ideas will be fruitful, leading to further analysis of searching and routing in graphs of all kinds.

**Acknowledgments.** Thanks to my advisers Olle Hẅggström and Devdatt Dubhashi, as well as Ian Clarke who originally suggested the edge updating algorithm, and Jon Kleinberg for taking the time to listen to and reflect upon my ideas.

DIVISION OF MATHEMATICAL STATISTICS
DEPARTMENT OF MATHEMATICAL SCIENCES
CHALMERS AND GÖTEBORG UNIVERSITY
412 96 GÖTEBORG
SWEDEN
E-MAIL: ossa@math.chalmers.se